\documentclass[11pt]{amsart}

\title{Biflatness of $\lp{1}$-semilattice algebras}
\author{Yemon Choi}
\date{10th May 2007}
\address{School of Mathematics and Statistics, University of Newcastle upon Tyne, Newcastle upon Tyne, NE1 7RU, England.}
\address{Current address: Department of Mathematics, University of Manitoba, Winnipeg, Manitoba R3T 2N2, Canada.}
\email{y.choi.97@cantab.net}
\thanks{\ackndiag}
\subjclass[2000]{Primary 46M20, 46J40; Secondary 43A20}
\usepackage{concrete}

\RequirePackage{amsmath,amsfonts,amsthm}

\RequirePackage[PS,nohug]{diagrams} 
\newarrow{Sub}       C--->  
\newarrow{Quot}    ----{>>} 

\newcommand{\ackndiag}{Uses Paul Taylor's {\tt diagrams.sty} macros.}

\newcommand{\YCrem}[1]{}   

\newcounter{pulse}
\numberwithin{pulse}{section}
\theoremstyle{plain}
\newtheorem{thm}[pulse]{Theorem}
\newtheorem{propn}[pulse]{Proposition}
\newtheorem{lemma}[pulse]{Lemma}
\newtheorem{coroll}[pulse]{Corollary}
\theoremstyle{definition}
\newtheorem{defn}[pulse]{Definition}
\newtheorem{eg}[pulse]{Example}
\newtheorem*{notn}{Notation}
\newtheorem{Rem}[pulse]{Remark}
\theoremstyle{remark}
\newtheorem*{rem}{Remark}

\renewcommand{\emph}[1]{{\sl #1\/}} 
\newcommand{\dt}[1]{{\it #1}}   
\newcommand{\defeq}{:=}
\newcommand{\st}{\,:\,}
\newcommand{\id}[1][]{{\sf 1}_{#1}} 

\newcommand{\abs}[1]{\vert{#1}\vert}
\newcommand{\Abs}[1]{\left\vert{#1}\right\vert}
\newcommand{\norm}[1]{\Vert{#1}\Vert}

\newcommand{\pair}[2]{\langle #1 ,\, #2\rangle}

\newcommand{\lp}[1]{\ell^{#1}}
\newcommand{\Lp}[1]{L^{#1}}
\newcommand{\tp}{\otimes}
\newcommand{\ptp}{\widehat{\otimes}}

\newcommand{\Cplx}{{\mathbb C}}
\newcommand{\Nat}{{\mathbb N}}
\newcommand{\Real}{{\mathbb R}}

\newcommand\Alg[1][]{{\mathcal A}_{#1}} 
\newcommand\Pt[1]{{\sf P}(#1)} 
\newcommand\Sch[1][]{{\sf  Sch}_{#1}}  
\newcommand\dset[2][]{(#2]_{#1}}  

\newcommand\Ind{{\mathbb I}}

\newcommand{\Schutzrep}{Sch{\"u}tzen\-berger representation}

\newcommand\MobF{M\"obius function}

\newcommand\MM{{\sf M}}  

\newcommand{\al}{\alpha}
\newcommand{\bt}{\beta}
\newcommand{\del}{\delta}
\newcommand{\kp}{\kappa}
\newcommand{\sig}{\sigma}
\newcommand{\tht}{\theta}
\newcommand{\om}{\omega}
\newcommand{\Om}{\Omega}
\newcommand\Und[1]{#1^\circ}

\newcommand{\bbG}{{\mathbb G}}  
\newcommand{\BG}[1][]{{\mathcal B}\/_{\bbG}}
\newcommand{\lpsum}[1]{\mathop{\overset{(\lp{#1})}{\bigoplus}}}

\begin{document}
\maketitle

\begin{section}{Introduction}
The notion of amenability for Banach algebras has received much attention since its introduction in the 1970s. If $G$ is a locally compact group, a celebrated theorem of B.~E.~Johnson states that the convolution algebra $\Lp{1}(G)$ is amenable if and only if $G$ is an amenable group.

$\lp{1}$-convolution algebras may also be defined for general (discrete) semigroups: here the problem of characterising amenability of the algebra in terms of the underlying semigroup seems not to admit such an elegant answer. One family of semigroups where a complete characterisation has been known for some time is the class of semilattices. If $S$ is a finite semilattice then its $\lp{1}$-convolution algebra $\Alg[S]$ is amenable. Conversely, if $S$ is an infinite semilattice, then by \cite[Theorem~10]{DuncNam} $\Alg[S]$ cannot be amenable.

A weaker notion than amenability is biflatness. In this paper we show that if $L$ is a semilattice then $\Alg[L]$ is biflat precisely when $L$ is ``uniformly locally finite'' (a notion that will be defined in due course). Our proof technique shows in passing that if $\Alg[L]$ is biflat then it is isomorphic as a Banach algebra to the Banach space $\lp{1}(L)$ equipped with \emph{pointwise} multiplication (this latter algebra is in fact known to be biprojective, a property of Banach algebras that is in general stronger than biflatness).

\begin{rem}
The original motivation for the work done here was somewhat indirect. In the article \cite{YC_GMJ} it was shown that for every semilattice $S$ the algebra $\Alg[S]$ is \dt{simplicially trivial}, that is, has vanishing simplicial cohomology in degrees $1$ and above. Since every biflat Banach algebra is  simplicially trivial (see e.g.~\cite[Propn 2.8.62]{Dal_BAAC}) it is natural in light of \cite{YC_GMJ} to enquire just when $\Alg[S]$ is biflat.
\end{rem}

\subsection*{Overview of this article}
We have tried to make this article reasonably accessible to both functional analysts and semigroup theorists. To keep our account to a reasonable length we have omitted several routine proofs which can instead be found in the references cited.

After some preliminary material, in Section \ref{s:diag-norm} we make some elementary observations about norms of diagonals for finite-dimensional algebras. In Section \ref{s:Schutz} we specialise to the case of the \Schutzrep\ of a locally finite semilattice, and in Section \ref{s:results} we assemble these observations to state and prove our main results (Theorems \ref{t:easydirection} and \ref{t:mainresult}).

We finish by \emph{sketching} how one can extend or build on the calculations of this article to characterise biflatness of the $\lp{1}$-convolution algebra of a Clifford semigroup in terms of its set of idempotents and its group components. The precise formulation is in Theorem~\ref{t:Cliff_biflat}. While this result generalises those for the semilattice case, we feel it is clearer and more natural -- at least in our approach -- to prove the special case first and then deduce the general case.
\end{section}

\begin{section}{Definitions and preliminaries}\label{s:prelim}

\subsection*{General notation}
We use $\abs{\quad}$ to denote the cardinality of a set. If $E$ is a (complex) vector space and $F$ a vector space of linear maps $E \to \Cplx$, we denote the corresponding pairing $F\times E \to \Cplx$ by $\pair{\quad}{\quad}$: thus given $x\in E$ and $\psi \in F$ we write $\pair{\psi}{x}$ rather than $\psi(x)$.

\subsection*{The partial order on a semilattice}
\begin{defn}
Recall that a semigroup $S$ is said to be a \dt{semilattice} if it is commutative and each element of $S$ is idempotent. The \dt{canonical partial order} on a semilattice $S$ is given by
\[ e\preceq f \iff ef=e \]
\end{defn}

We recall for later convenience that the \dt{greatest lower bound} or \dt{meet} of two elements $s, t \in S$, with respect to this partial order, is nothing but their product $st$. In particular, a minimal element of the partially ordered set $(S,\preceq)$ is nothing but a \dt{zero element} of $S$, i.e.~ an element $\tht$ satisfying $x\tht=\tht=\tht x$ for all $x \in S$.

\begin{notn}
If $(P,\preceq)$ is a partially ordered set we let
\[ \dset[P]{x} \defeq \{ y \in P \st y \preceq x \} \]
We shall sometimes drop the subscript and write $\dset{x}$ if it is clear which partially ordered set we are referring to.
\end{notn}

\begin{defn}\label{dfn:downset}
Let $P$ be a partially ordered set and let $D \subseteq P$. We say that $D$ is a \dt{downwards-directed set} or \dt{downset} in $P$ if it has the following property:
\[ \text{ if $x \in D$, $y \in P$ and $y\preceq x$ then $y \in D$} \]
\end{defn}

The following observations are immediate from the definition of the canonical partial order.
\begin{lemma}\label{l:downset-ideal}
Let $S$ be a semilattice, equipped with the canonical partial order. If $x, y\in S$ then $x\preceq y$ if and only if $y$ divides $x$, i.e.~if and only if there exists $z \in S$ with $x=yz$. Hence for each $f \in S$, $\dset[S]{f} = fS=Sf$ is just the ideal generated by $f$.

More generally, any downset in $S$ is an ideal in $S$, hence is in particular a subsemigroup of $S$.
\end{lemma}

We are interested in how the structure of a semilattice $S$ is reflected in properties of its \dt{$\lp{1}$-convolution algebra}, which will be denoted throughout by $\Alg[S]$. (Recall that $\Alg[S]$ is the Banach algebra obtained by completing the semigroup algebra $\Cplx S$ with respect to the $\lp{1}$-norm; see \cite[Example 1.23]{Bons-Dunc} or \cite[Example 2.1.23(v)]{Dal_BAAC} for further details).

\begin{notn}
 The canonical basis vectors in $\Alg[S]$ will be denoted by $(e_s)_{s\in S}$, where $e_s$ is the function $S\to\Cplx$ that is $1$ at $s$ and $0$ everywhere else.

By definition of the convolution product in $\Alg[S]$, $e_se_t=e_{st}$ for all $s,t \in S$.
\end{notn}

\subsection*{Biflatness and amenability for Banach algebras}
Most of the background on Banach algebras, modules and tensor products that is needed may be found in \cite{Bons-Dunc}: in particular if $X$ and $Y$ are Banach spaces we denote their \dt{projective tensor product} by $X\ptp Y$. We differ from the notation in \cite{Bons-Dunc} slightly in that we denote the dual of a Banach space $E$ by $E'$; similarly if $f:E\to F$ is a bounded linear map between Banach spaces then the \dt{adjoint} map $F'\to E'$ will be denoted by $f'$. (This follows the notation in \cite{Dal_BAAC} for instance.)

To avoid cumbersome repetition, we adopt the convention that if $X$, $Y$ are Banach $A$-modules, an \dt{$A$-module map from $X$ to $Y$} will \emph{always} mean a bounded linear map $X\to Y$ which respects the $A$-action.

\begin{notn}
If $A$ is a Banach algebra then the linearised product map $\pi:A\ptp A\to A$ is defined by
\[ \pi(a\tp b) \defeq ab\qquad(a,b\in A) \]
\end{notn}

Although we shall not need to consider amenability in any detail, we shall in proving our main result make use of the notion of a \dt{virtual diagonal}.
\begin{defn}[see {\cite[\S43]{Bons-Dunc}}]\label{dfn:vd}
Let $A$ be a Banach algebra: then $A$ and $A\ptp A$, and hence their respective second duals, may be regarded as $A$-bimodules in a canonical way. A \dt{virtual diagonal} for $A$ is an element $M\in(A\ptp A)''$ with the following properties:
\begin{itemize}
\item $aM =Ma$
\item $a\pi''(M)=a=\pi''(M)a$
\end{itemize}
for all $a \in A$.

A Banach algebra with a virtual diagonal is said to be \dt{amenable}.
\end{defn}

The notion of biflatness is perhaps less well-known than that of amenability. There are several equivalent characterisations of biflatness for Banach algebras; we shall use the following one.
\begin{defn}[cf.\ {\cite[Exercise VII.2.8]{Hel_HBTA}}]
\label{dfn:biflat}
Let $A$ be a Banach algebra. $A$ is \dt{biflat} if there exists an $A$-bimodule map $\sig: A\to (A\ptp A)''$ such that the diagram
\[ \begin{diagram}[tight,height=2.5em]
A & \rTo^\sig & (A\ptp A)'' \\
 & \rdTo_{\kp} & \dTo_{\pi''} \\
 & & A''
\end{diagram} \]
commutes, where $\kp: A \to A''$ is the natural embedding of $A$ in its second dual.

Given $C \geq 1$, we say that $A$ is \dt{$C$-biflat} if the map $\sig$ can be chosen to have norm less than or equal to $C$.
\end{defn}

\begin{Rem}\label{r:VD_from_1}
If $A$ is biflat and unital, with identity $\id$, then it is easily checked that $\sig(\id)$ is a virtual diagonal for $A$; in particular, $A$ is amenable.
\end{Rem}

We shall not discuss the precise connection between amenability and biflatness here. However, the following example will play a key role in what follows.

\begin{eg}[Example of a non-amenable, biflat Banach algebra]\label{eg:ptwise_l1}
Let $\Om$ be any set and let $(\del_\om)_{\om\in\Om}$ denote the standard basis of $\lp{1}(\Om)$. We may equip $\lp{1}(\Om)$ with pointwise multiplication, i.e.
\[ \del_\al \del_\bt \defeq \left\{ \begin{aligned}
	\del_\al & \quad\text{ if $\al=\bt$}\\
	0 & \quad\text{ otherwise}
\end{aligned}\right. \]
This makes $\lp{1}(\Om)$ into a Banach algebra, which will be denoted throughout by $\Pt{\Om}$.

\emph{We claim that this Banach algebra is biflat}: indeed, the map $\sig: \Pt{\Om}\to (\Pt{\Om}\ptp\Pt{\Om}) ''$ defined by
\[ \sig(\del_\om) \defeq \del_\om \tp \del_\om \quad\quad(\om\in\Om) \]
satisfies the properties required in Definition \ref{dfn:biflat}. On the other hand, if $\Om$ is infinite it is well-known that $\Pt{\Om}$ cannot be amenable
(as
 it has no bounded approximate identity; see \cite[Theorem~2.9.65 and Example~4.1.42]{Dal_BAAC} for details).
\end{eg}

\begin{rem}
The argument above shows that $\Pt{\Om}$ is in fact a \dt{biprojective} Banach algebra: that is, the map $\sig$ in the definition of biflatness can be chosen to take values in $\Pt{\Om}\ptp\Pt{\Om}\subseteq (\Pt{\Om}\ptp \Pt{\Om})''$. We shall not discuss biprojectivity in this article.
\end{rem}

\subsection*{Properties of $\lp{1}$}
It is evident from the definitions above that a concrete picture of the $\Alg[S]$-bimodule  $\Alg[S]\ptp\Alg[S]$ will be helpful for our investigations. Since the underlying Banach space of $\Alg[S]$ is $\lp{1}$, we pause to recall some trivial but useful observations about $\lp{1}$-spaces.

\begin{Rem}\label{r:bidual-of-incl}
If $X$ is a non-empty subset of $Y$ then we may identify $\lp{1}(X)$ with a closed, complemented subspace of $\lp{1}(Y)$, namely the space of all elements of $\lp{1}(Y)$ whose support is contained in $X$. We may go on to identify the bidual $\lp{1}(X)''$ with a closed complemented subspace of $\lp{1}(Y)''$, namely the \emph{annihilator} of $\lp{\infty}(Y\setminus X)$
 viewed as a subspace of $\lp{\infty}(Y)$.
\end{Rem}

If $X$ and $Y$ are two sets then it is well known that there is an isometric linear isomorphism $\lp{1}(X)\ptp\lp{1}(Y) \to \lp{1}(X\times Y)$: see e.g.~\cite[Examples~2.1.24]{Dal_BAAC}.
Thus the underlying Banach space of $\Alg[S]\ptp\Alg[S]$ is $\lp{1}(S\times S)$, and
 the bimodule action on $\lp{1}(S\times S)$ is uniquely defined by requiring that
\[ e_x\cdot e_{(s,t)} \defeq e_{(xs,t)} \quad\text{ and }\quad
	e_{(s,t)}\cdot e_x \defeq e_{(s,tx)} \qquad(s,t,x\in S) \]
and extending by linearity and continuity.

In what follows we shall repeatedly use the following well-known observation, without further comment:

\smallskip
\noindent{\bf Fact.} Let $E$ be a Banach space, let $\Om$ be a set and let $T:\lp{1}(\Om)\to E$ be a bounded linear map. Then
\[ \norm{T}=\sup_{\om\in\Om} \norm{T(\del_\om)} \]
where $(\del_\om)_{\om\in\Om}$ is the canonical unit basis of $\lp{1}(\Om)$. Conversely, if $f:\Om\to E$ is a bounded function, there exists a unique bounded linear map $\widetilde{f}:\lp{1}(\Om)\to E$ which satisfies
\[ \widetilde{f}(\del_\om)=f(\om) \qquad\text{ for all $\om\in\Om$.} \]

\smallskip
The proof of this observation needs only the triangle inequality and the definition of the $\lp{1}$-norm, and is left to the reader.

\end{section}

\begin{section}{Isomorphism constants and diagonals in the finite dimensional setting}
\label{s:diag-norm}
\begin{defn}
Let $A$ be a complex algebra, and let $\pi:A\tp A\to A$ denote the linearisation of the product map. An element $\Delta\in A\tp A$ is called a \dt{diagonal} for $A$ if
\begin{itemize}
\item $a\cdot \Delta = \Delta\cdot a$
\item $\pi(\Delta)a=a=a\pi(\Delta)$
\end{itemize}
for all $a \in A$.
\end{defn}

\begin{rem}
If $A$ is a \emph{finite-dimensional} amenable Banach algebra, then any virtual diagonal for $A$ is in fact a diagonal for $A$ (since $(A\ptp A)''=A\tp A$).
\end{rem}

Let $\Phi$ be any finite set and equip the vector space $\Cplx^\Phi$ with pointwise multiplication. Then $\Cplx^\Phi$ has a \emph{unique}\footnotemark\ diagonal, namely
$\Delta_\Phi= \sum_{i \in \Phi} \del_i\tp \del_i$
where $\del_x$ denotes the standard basis vector in $\Cplx^\Phi$ corresponding to the element $x \in \Phi$.
\footnotetext{Uniqueness is easily checked by solving the identities $\delta_j M=M\delta_j$ for $j\in\Phi$ and $\pi(M)=\sum_{i\in\Phi}\delta_i$, where $M$ is a $\Phi\times\Phi$ matrix. In this context $\pi$ corresponds to `restriction of a matrix to its diagonal entries'.}

Suppose furthermore that $A$ is another complex algebra and that we have an algebra isomorphism $\al: \Cplx^\Phi \to A$. Then $A$ has a unique diagonal $\Delta\defeq\sum_{x \in \Phi} \al(\del_x)\tp\al(\del_x)$.

Pick a basis for $A$, say $\{ e_i : i \in \Ind \}$ for some index set $\Ind$. Suppose also that with respect to the bases $\{\del_x\}$ and $\{e_i\}$ the linear map $\al$ has a real-valued matrix, i.e.
\[ \al(\del_x) = \sum_{s \in \Ind} M(s,x)e_s \]
for some function $M: \Ind\times\Ind \to \Real$. Then
\[ \begin{aligned}
\Delta & =  \sum_{x \in \Phi} \al(\del_x)\tp\al(\del_x)  \\
 & = \sum_{x \in \Phi} \left(\sum_{s\in \Ind } M(s, x)e_s\right)\tp\left(\sum_{t\in \Ind} M(t,x)e_t\right) \\
 & = \sum_{(s,t) \in \Ind\times \Ind} \left(\sum_{x \in\Phi} M(s,x)M(t,x) \right) e_s\tp e_t
\end{aligned} \]

We apply this as follows. Equip $A$ with the $\lp{1}$-norm with respect to the basis $(e_i)_{i\in\Ind}$ (this need not be an \emph{algebra norm} on $A$, although in the cases of interest to us it will be). We then equip $A\tp A$ with the corresponding $\lp{1}$-norm (with respect to the basis $(e_s\tp e_t)_{s,t\in\Ind}$).

There is then a crude lower bound
\begin{equation}\label{eq:lower_bound_on_diag}
\begin{aligned}
\norm{\Delta}_1 & = \sum_{(s,t) \in \Ind\times \Ind} \Abs{\sum_{x \in \Phi} M(s,x) M(t,x) } &  \\
 & \geq \sum_{s \in \Ind} \Abs{\sum_{x\in\Phi} M(s,x)M(s,x)} 
 & = \sum_{s \in \Ind} \sum_{x \in\Phi} M(s,x)^2
\end{aligned}
\end{equation}
Therefore, if we equip $\Cplx^\Phi$ with the $\lp{1}$-norm (with respect to the basis $(\del_x)$) we find using the Cauchy-Schwarz inequality that
\begin{equation}\label{eq:bound1}
 \begin{aligned}
\norm{\al} 
	 = \sup_{x\in\Phi} \norm{\al(\del_x)}_1 
	& = \sup_{x\in\Phi} \sum_{s\in\Ind} \abs{M(s,x)} \\
	& \leq \sup_{x\in\Phi} \abs{\Ind}^{1/2} \left(\sum_{s\in\Ind}M(s,x)^2\right)^{1/2} \\
	& \leq \abs{\Ind}^{1/2} \left( \sum_{x\in\Phi}\sum_{s\in\Ind}M(s,x)^2 \right)^{1/2} \\
	& \leq \abs{\Ind}^{1/2}\norm{\Delta}^{1/2}
\end{aligned} \end{equation}

We shall apply these general estimates below to the special case of the \dt{\Schutzrep} of a finite semilattice. The full definition will be given in the next section: we remark that in this case one has an algebra isomorphism $\Sch: \Alg[L] \to \Cplx^L$, and if in the above we take $\Phi=L$, $A=\Alg[L]$ and $\al=\Sch^{-1}$ then $M$ turns out to be the \dt{\MobF} for $L$.
\end{section}

\begin{section}{The \Schutzrep}\label{s:Schutz}
The calculations that follow are inspired by the presentation in \cite{St_Mob}, although our notation will differ slightly. We remark that the constructions in \cite{St_Mob} are much more general than is needed for our work: most of what we need has been known for over 30 years and goes back to results of Solomon (see also \cite{Gr_MA} for a concise account).

\begin{notn}
If $S$ is a semilattice we shall write $\Und{S}$ for the underlying set of $S$.
\end{notn}

\begin{defn}
Let $S$ be a semilattice equipped with its canonical partial ordering.
 We define a bounded linear map $\Sch[S]: \Alg[S] \to \lp{\infty}(\Und{S})$ by
\[ \Sch[S](e_t) \defeq {\mathbf 1}_{\dset[S]{t}} \quad\quad(t\in S) \]
where ${\mathbf 1}_X$ denotes the indicator function of a subset $X\subseteq \Und{S}$.

$\Sch[S]$ is the \dt{\Schutzrep} of $S$. If it is clear which semilattice we are working with we shall occasionally suppress the $S$ and simply write $\Sch$.
\end{defn}
Crucially for what follows $\Sch[S]$ is an algebra homomorphism (the easiest way to see this is to check that $\Sch[S](e_{st})=\Sch[S](e_s)\Sch[S](e_t)$ for all $s,t\in S$). See, e.g.~\cite{Gr_MA} for more on the properties of $\Sch$ (at least for the case where $S$ is a \emph{finite} semilattice).

\begin{rem}
In effect $\Sch$ represents the convolution algebra $\Alg[S]$ as an algebra of functions on some carrier space. In this context, note that the Gelfand transform of $\Alg[S]$ has the form $\Alg[S] \to C_0(\widehat{S})$
where $\widehat{S}$ is the space of characters on $S$. Now for any semilattice there is a canonical inclusion $\Und{S} \rSub \widehat{S}$,
giving us a ``restriction map''  $\lp{\infty}(\widehat{S}) \to \lp{\infty}(\Und{S})$; the composition
\[ \Alg[S] \rTo C_0(\widehat{S})\rSub \lp{\infty}(\widehat{S})
   \rTo \lp{\infty}(\Und{S}) \]
is just $\Sch$.
 \end{rem}

\begin{notn}
We denote basis vectors in $c_0(\Und{S})$ by $\del_x$, for $x \in S$: that is, $\del_x(t)$ is defined to be $1$ when $t=x$ and
 $0$ otherwise.
\end{notn}
Note that if $L$ is a \emph{finite} unital semilattice then
\begin{equation}\label{eq:finitesum}
\Sch[L](e_t) = \sum_{s \in \dset[L]{t}} \del_s \quad\quad(t \in L)
\end{equation}
and we may regard $\Sch[L]$ an algebra homomorphism $\Alg[L]\to\Cplx^{L}$. A direct computation then shows that $\Sch[L]$ is invertible, with
\begin{equation}\label{eq:invSch}
\Sch[L]^{-1}(\del_x) = \sum_{y \in\dset[L]{x}} \mu(y,x) e_y
\end{equation}
where $\mu$ is the \dt{\MobF} for $L$. (For the definition and basic properties of the \MobF\ of a locally finite poset, see \cite[\S\S3.6--3.7]{Stan_EC1}.)

\subsection*{Norm estimates for $\Sch:\Alg[L]\to\Pt{\Und{L}}$}
Note that since $L$ is a finite semilattice it has a least element, which we denote by $\tht$. We define a family of idempotents $(\rho_t)_{t\in L}$ by
\[ \rho_t = \left\{ \begin{aligned}
	e_{\tht} & \quad\text{ if $t=\tht$} \\
	\prod_{x\in L \st x \prec t} (e_t-e_x) & \quad\text{ otherwise.}
 \end{aligned} \right. \]

\begin{propn}[{cf.~\cite[Propn.~1.5]{Eti_Mob}}]\label{p:magic_idempotents}
$\Sch(\rho_t) = \del_t$.
\end{propn}

The proof given in \cite{Eti_Mob} of this identity leaves some of the details to the reader, and so for sake of completeness we give a full proof in the Appendix. In any case, all we need is the following easy consequence.

\begin{coroll}\label{c:crudebound}
$\norm{\Sch^{-1}}\leq 2^{\abs{L}-1}$.
\end{coroll}

\begin{proof}[Proof of corollary]
Let $t \in L$: then by Proposition \ref{p:magic_idempotents},
\[ \norm{\Sch^{-1}(\del_t)} = \norm{\rho_t}
 \leq \prod_{x\in L \st x \prec t} \norm{e_t-e_x}
 =2^{\Abs{\dset[L]{x}-1}} \leq 2^{\abs{L}-1}\;. \]
Since $t$ is arbitrary, $\max_t \norm{\Sch^{-1}(\del_t)}\leq 2^{\abs{L}-1}$ and the result follows.
\end{proof}

Thus $\norm{\Sch}$ and $\norm{\Sch^{-1}}$ are controlled by the size of our finite semilattice $L$. Later on we shall need to know that they are also controlled by the norm of any diagonal in $\Alg[L]$; this is made precise by the following result.

\begin{propn}\label{p:control-via-diag}
Let $L$ be a finite unital semilattice, and suppose $\Alg[L]$ has a diagonal with norm $\leq C$. Regard the \Schutzrep\ $\Sch$ as a linear map from the normed space $\Alg[L]$ to the normed space $\Pt{\Und{L}}$.
Then
\[ \norm{\Sch}=\abs{L}\leq C \quad\text{ and } \quad \norm{\Sch^{-1}}\leq C \;.\]
\end{propn}
\begin{proof}
We shall apply the estimates from the previous section with $A=\Alg[L]$ and $\al=\Sch^{-1}$.
Since $L$ is a finite semilattice, it is immediate from Equation \eqref{eq:finitesum} that
\[ \norm{\Sch} =\sup_{t \in L} \norm{\Sch(e_t)}_1 = \sup_{t \in L} \abs{ \dset{t} } = \abs{\dset{\id[L]}}= \abs{L} \;;\]
and by the estimate \eqref{eq:bound1}, $\norm{\Sch^{-1}}\leq\abs{L}^{1/2}C^{1/2}$. Therefore it remains only to show that $\abs{L}\leq C$. But from the original bound \eqref{eq:lower_bound_on_diag} we know that
\[ C \geq \norm{\Delta} \geq \sum_{s \preceq x}\mu(s,x)^2 \geq \sum_{x\in L} \mu(x,x)^2 \]
Since the \MobF\ $\mu$ satisfies $\mu(t,t)=1$ for all $t \in L$ this concludes the proof.
\end{proof}
\end{section}

\begin{section}{Results}\label{s:results}
We recall the following result of Duncan and Namioka, 
which will be needed later.
\begin{thm}[{\cite[Theorem 10]{DuncNam}}]\label{t:DuncNam}
Let $L$ be a semilattice. Suppose that $\Alg[L]$ is amenable, i.e.~possesses a virtual diagonal. Then $L$ is finite.
\end{thm}

A partially ordered set $(P,\preceq)$ is said to be \dt{locally finite} if $\dset{x}$ is finite for each $x \in P$. It is natural to then consider the following stronger notion.
\begin{defn}
We say that a semilattice $L$ is \dt{locally $C$-finite}, for some constant $C>0$, if  $\abs{\dset{f}} \leq C$ for all $f \in S$. A semilattice which is locally $C$-finite for some $C$ is said to be \dt{uniformly locally finite.}
\end{defn}

\begin{eg}
Some illustrative examples:
\begin{enumerate}
\item Every finite semilattice is uniformly locally finite.
\item Let $E$ be a semilattice consisting of a zero element $\theta$ and infinitely many orthogonal idempotents, i.e.~$xy=\theta=yx$ for all distinct $x,y\in E\setminus\{\theta\}$. Then $E$ is locally $2$-finite.
\item All uniformly locally finite semilattices have finite height, where the \dt{height} of a partially ordered set is simply the supremum over all lengths of chains in that set. The converse is false: take the semilattice $E$ from the previous example and adjoin an identity element to get a semilattice of height 2 which is not even locally finite.
\item Let $(\Nat, \min)$ denote the semilattice with underlying set $\Nat$ and multiplication given by $(m,n)\mapsto \min(m,n)$. Then the partial order on $(\Nat,\min)$ is just the usual order on $\Nat$; in particular $(\Nat, \min)$ is locally finite but not uniformly locally finite.
\end{enumerate}
\end{eg}

In view of Example (2) above, one might hope that a uniformly locally finite semilattice is close to being finite, apart from a possibly infinite set of maximal elements. The following example shows that things can be a little more complicated.

\begin{eg}
Let $S$ be the semilattice defined as follows: the underlying set is $(\Nat\times\{1,2\})\sqcup\{\tht\}$, while multiplication is defined by taking $\theta$ to be a zero element and setting
\[ (m,i)\cdot (n,j) \defeq 
\left\{\begin{aligned} (r,2) & \quad\text{ if $m=2r-1$, $n=2r$ and $i=j=1$} \\ \tht & \quad\text{otherwise.} \end{aligned}\right.  \]
The maximal elements of $S$ are all those of the form $(m,1)$ for $m \in \Nat$; in particular $S$ contains infinitely many non-maximal elements. On the other hand a quick calculation shows that $S$ is locally $2$-finite.
\end{eg}

Despite this variety of possible uniformly locally finite semilattices, it turns out that the $\lp{1}$-convolution algebra of such a semilattice is determined up to isomorphism by the cardinality of the semilattice! This is a consequence of the first of our main results, which is as follows.

\begin{thm}\label{t:easydirection}
Let $S$ be a locally $C$-finite semilattice, for some $C\geq1$, and let $\Sch:\Alg[S]\to\lp{\infty}(\Und{S})$ denote the \Schutzrep\ of $S$. Then:
\begin{enumerate}
\item[(i)] $\norm{\Sch(a) }_1 \leq C \norm{a}$ for all $a \in \Alg[S]$;
\item[(ii)] $\Sch: \Alg[S] \to \Pt{\Und{S}}$ is invertible, and
\[ \norm{\Sch^{-1}(b)} \leq 2^{C-1}\norm{b}_1 \quad\quad\text{ for all $b \in \Pt{\Und{S}}$} \;. \]
\end{enumerate}
In particular, $\Alg[S]$ is biflat (since $\Pt{\Und{S}}$ is).
\end{thm}

We thus obtain plenty of examples of semilattices $S$ for which $\Alg[S]$ is biflat. The thrust of our second main result is that these are \emph{all} the possible examples.

\begin{thm}\label{t:mainresult}
Let $S$ be a semilattice such that $\Alg[S]$ is $C$-biflat for some constant $C\geq1$. Then $S$ is locally $C$-finite. 

Also, if we regard the \Schutzrep\ of $S$ as an algebra homomorphism $\Sch: \Alg[S] \to \Pt{\Und{S}}$, then $\norm{\Sch^{-1}}\leq C$.
\end{thm}

Before proving these results we record some lemmas.
\begin{lemma}\label{l:cut-down}
Fix $x \in S$. Then
 $e_x\Xi e_x \in \left(\Alg[xS]\ptp\Alg[xS]\right)''$ for all $\Xi\in \left(\Alg[S]\ptp\Alg[S]\right)''$.
\end{lemma}
\begin{proof}\
Let $T=S\setminus xS=S\setminus Sx$.
By Remark \ref{r:bidual-of-incl} we can identify 
 $\lp{1}(xS\times Sx)''$ with the annihilator (in $\lp{\infty}(S\times S)'$) of $\lp{\infty}(T\times T)$.
Therefore it suffices to show that $\pair{e_x\Xi e_x}{\psi}=0$ for all $\psi \in \lp{\infty}(T\times T)$.

Since $\pair{e_x\Xi e_x}{\psi}\equiv \pair{\Xi}{e_x\psi e_x}$
it suffices to show that $e_x\psi e_x =0$ for all $\psi \in \lp{\infty}(T\times T)$. But by definition of the $\Alg[S]$-action on $\lp{\infty}(S\times S)$, for every $s \in S$ we have
\[ (e_x\psi e_x)_s = \psi_{xsx} = 0 \]
since $xsx \in xS=S\setminus T$ and $\psi$ is supported on $T$.
\end{proof}

\begin{lemma}\label{l:Schutz-restr}
Let $S$ be a semilattice and let $\Sch[S]$ be its \Schutzrep. If $L\subseteq S$ is a downset in $S$, then $\Sch[S]\left(\Alg[L]\right)\subseteq \lp{\infty}(\Und{L})$
where we view $\lp{\infty}(\Und{L})$ as a subspace of $\lp{\infty}(\Und{S})$.
Moreover
\[\Sch[L]=\Sch[S]\vert_{\Alg[L]}\;.\]
\end{lemma}
\begin{proof}
This is immediate from the definition of the \Schutzrep\
and Lemma \ref{l:downset-ideal}.
\end{proof}

\begin{proof}[Proof of Theorem \ref{t:easydirection}]
Suppose that $S$ is locally $C$-finite.
If $a = \sum_{f \in S} a_f e_f \in \Alg[S]$ then
\[ \norm{\Sch(a)}_1 \leq \sum_{f \in S} \abs{a_f} \norm{\Sch(e_f)}_1
= \sum_{f \in S} \abs{a_f} \abs{\dset{f}} \leq C\norm{a} \;; \]
hence $\Sch$ takes values in $\lp{1}(\Und{S})$, and has norm $\leq C$ when regarded as a linear map from $\Alg[S]$ to $\lp{1}(\Und{S})$.

Thus $\Sch:\Alg[S]\to\Pt{\Und{S}}$ is a well-defined, continuous algebra homomorphism.

We shall now construct a 
bounded linear
map $\MM : \lp{1}(\Und{S}) \to \Alg[S]$
which has
 norm $\leq 2^{C-1}$ and which is a $2$-sided inverse to $\Sch$. If $x \in \Und{S}$ let
\[ \MM(\del_x) \defeq (\Sch[\dset{x}])^{-1} \del_x \]

By Corollary \ref{c:crudebound},
\[ \norm{\MM(\del_x)} \leq 2^{\abs{\dset{x}}-1} \leq 2^{C-1} \quad\text{ for all $x\in\Und{S}$;}\]
thus $\MM$ extends to a bounded linear map $\lp{1}(\Und{S}) \to \Alg[S]$ which has norm $\leq 2^{C-1}$. For any $x \in \Und{S}$ we have (using Lemma \ref{l:Schutz-restr})
\[ \Sch\MM(\del_x) = \Sch (\Sch[\dset{x}])^{-1} \del_x = \Sch[\dset{x}](\Sch[\dset{x}])^{-1} \del_x = \del_x \quad; \]
thus by linearity and continuity $\MM$ is right inverse to $\Sch$.

It remains to show that $\MM$ is left inverse to $\Sch$. Let $f \in S$: then for any $x$ such that $x\preceq f$ we have $\dset{x} \subseteq \dset{f}$, and so by Lemma \ref{l:Schutz-restr} the restriction of $\Sch[\dset{f}]^{-1}$ to $\Pt{\Und{\dset{x}}}$ is just $\Sch[\dset{x}]^{-1}$. Hence
 (appealing again to Lemma \ref{l:Schutz-restr} ) we have
\[ \begin{aligned}
\MM\Sch(e_f) = \MM\left( \sum_{x\in\dset{f}} \del_x \right) 
& = \sum_{x\in\dset{f}} (\Sch[\dset{x}])^{-1}(\del_x) \\
& = \sum_{x\in\dset{f}} (\Sch[\dset{f}])^{-1}(\del_x) \\
& =  (\Sch[\dset{f}])^{-1} \left( \sum_{x\in\dset{f}} \del_x \right) \\
	& =  (\Sch[\dset{f}])^{-1} \Sch[\dset{f}](e_f) & = e_f\;.
\end{aligned} \]
By linearity and continuity we deduce that $\MM\Sch(a)=a$ for each $a \in \Alg[S]$,  and this concludes the proof.
\end{proof}
\begin{rem}
We could have shortened the proof slightly if we
had known
 beforehand that $\Sch$ is injective (since then it suffices to show that the map $\MM$ is a right inverse to $\Sch$, without explicitly showing it is also a left inverse). It would be interesting to know for which semilattices the \Schutzrep\ is injective.
\end{rem}

\begin{proof}[Proof of Theorem \ref{t:mainresult}]
Suppose $\Alg[S]$ is $C$-biflat: then there exists an $\Alg[S]$-bimodule map $\sig: \Alg[S] \to \left(\Alg[S]\ptp\Alg[S]\right)''$ such that $\pi''\sig=\kp$ and $\norm{\sig}\leq C$.

Fix $f \in S$. Then for all $x\in S$,
\[  \sig(e_{fx}) = \sig(e_fe_xe_f) =e_f\sig(e_x)e_f \]
so by Lemma \ref{l:cut-down} $\sig(e_{fx})\in\left(\Alg[\dset{f}]\ptp\Alg[\dset{f}]\right)''$. Recall that $\dset{f}$ is a unital semilattice and a subsemigroup of $S$. Hence the restriction $\sig_f$ of $\sig$ to $\Alg[\dset{f}]$ is an $\Alg[\dset{f}]$-bimodule map with values in $(\Alg[\dset{f}]\ptp\Alg[\dset{f}])''$, which makes the following diagram commute:
 \begin{diagram}[tight,height=2.5em]
\Alg[\dset{f}] & & \rTo^{\sig_f} & \left(\Alg[\dset{f}]\ptp\Alg[\dset{f}]\right)'' \\
 & \rdTo(3,2)_{\kp} & & \dTo_{\pi''} \\
 & & & \Alg[\dset{f}]''
\end{diagram} 

Since $e_f$ is an identity element for $\Alg[\dset{f}]$, $\sig_f(e_f)$ is a virtual diagonal for $\Alg[\dset{f}]$ (see Remark~\ref{r:VD_from_1}). By Theorem \ref{t:DuncNam} we deduce that $\dset{f}$ is finite, and so $\sig_f(e_f)$ is a diagonal for $\Alg[\dset{f}]$. As $\norm{\sig_f(e_f)}\leq \norm{\sig}\leq C$, we see from Proposition \ref{p:control-via-diag} that $\abs{\dset{f}}\leq C$. This holds for all $f\in S$ and thus $S$ is locally $C$-finite.

As in (the proof of) Theorem \ref{t:easydirection} we deduce that $\Sch$ is an algebra isomorphism from $\Alg[S]$ \emph{onto} $\Pt{\Und{S}}$; to finish it is enough to show that the map $\MM$ defined in the proof of Theorem \ref{t:easydirection} 
in fact has norm $\leq C$. But since by definition
\[ \MM(\del_x) \defeq (\Sch[\dset{x}])^{-1} \del_x \qquad(x\in \Und{S}) \;,\]
Proposition \ref{p:control-via-diag} implies that $\norm{\MM(\del_x)} \leq C$ for all $x\in\Und{S}$, and hence $\norm{\MM}\leq C$ as required.
\end{proof}

\end{section}

\section{A sketch of an extension to Clifford semigroups}\label{s:Clifford_sketch}
We have seen that methods from the representation theory of (locally) finite semilattices are useful for studying the $\lp{1}$-convolution algebras of such semigroups. One natural next step is to use the more general machinery for representations of finite semigroups to attack the question of biflatness for $\lp{1}(S)$ for various classes of $S$, in particular when $S$ is an \emph{inverse semigroup} (cf.~the techniques in \cite{St_Mob}).
Originally it was planned to address this in future work; the author has since been informed by F. Habibian \cite{FHab_email} that he is working on the same problem from a different approach, and so the author's preliminary calculations will not be submitted for publication. 

It nevertheless seems worthwhile to sketch how the methods of this paper adapt to the case of \dt{Clifford semigroups}. We therefore conclude this paper with a rapid account and direct the reader to the aforementioned work of Habibian for results on more general inverse semigroups. Details will be omitted in several places.

\subsection*{Notation, preliminaries}
Let $\bbG$ be a Clifford semigroup, $L$ its set of idempotents, and $\bbG=\coprod_{x \in L} G(x)$ the associated decomposition of $\bbG$ as a strong semilattice of groups (see e.g.~\cite[Ch.~4]{How_fund-sgp} for basic definitions and properties).

A basic fact which we shall use is that for each $x\in L$, $x$ lies in $G(x)$ and is the identity element for $G(x)$. We also need the fact that if $e,f\in L$ then $G_e\cdot G_f\subseteq G_{ef}$; in particular the function $q:\bbG \to L$ defined by
\begin{equation}\label{eq:retract_onto_L}
q(G_e) =\{e\} \qquad\text{ for all $x\in L$}
\end{equation}
is a semigroup homomorphism. Note that the restriction of $q$ to the embedded copy of $L$ inside $\bbG$ is just the identity map, and that
\[ tq(t) =t =q(t)t \qquad\text{ for all $t\in\bbG$.} \]

Our promised characterisation of the biflatness of $\Alg[\bbG]$ is as follows.
\begin{thm}\label{t:Cliff_biflat}
$\Alg[\bbG]$ is biflat if and only if both the following conditions hold:
\begin{itemize}
\item[$(i)$] $L$ is uniformly locally finite;
\item[$(ii)$] $G(x)$ is amenable for each $x\in L$.
\end{itemize}
\end{thm}
Our strategy for proving Theorem \ref{t:Cliff_biflat} hinges on two main propositions (\ref{p:Schutz-is-iso} and \ref{p:prelim_Cliff_biflat} below).
To state these results we introduce another Banach algebra built out of $\bbG$, namely the $\lp{1}$-sum
\[ \BG\defeq \lpsum{1}_{x\in \Und{L}} \Alg[G(x)] \]
in which multiplication is defined `componentwise' and multiplication in the $x$th component is just given by the usual multiplication in $\Alg[G(x)]$.

\begin{propn}[cf.\ Theorem \ref{t:easydirection}]\label{p:Schutz-is-iso}
Suppose $L$ is uniformly locally finite. Then there is a Banach algebra isomorphism from $\Alg[\bbG]$ onto $\BG$.
\end{propn}

\begin{propn}[cf.\ Example \ref{eg:ptwise_l1}]\label{p:prelim_Cliff_biflat}
$\BG$ is biflat if and only if $G(x)$ is amenable for all $x\in L$. 
\end{propn}

Before sketching the proofs of these propositions we show how they combine to yield Theorem \ref{t:Cliff_biflat}.

\begin{defn}
Given a Banach algebra $B$ and a closed subalgebra $A$, we say that $A$ is a \dt{retract} of $B$ if there exists a continuous homomorphism $B\to A$ which restricts to the identity on $A$.
\end{defn}

One can show by straightforward diagram-chasing that a retract of a biflat Banach algebra is also biflat (we omit the details).

\begin{proof}[Proof of Theorem \ref{t:Cliff_biflat}, assuming Pro\-pos\-itions \ref{p:Schutz-is-iso} and \ref{p:prelim_Cliff_biflat}]
Suppose that $(i)$ and $(ii)$ hold. By $(ii)$ the Banach algebra $\BG$ is biflat; and by $(i)$ $\Alg[\bbG]$ is isomorphic as a Banach algebra to $\BG$, so is itself biflat.

Conversely, suppose $\Alg[\bbG]$ is biflat. We \emph{claim} that $\Alg[L]$ is also biflat. The simplest way to see this is to note that it is a retract of $\Alg[\bbG]$: a suitable homomorphism $\Alg[\bbG]\to\Alg[L]$ can be defined by $e_t\mapsto e_{q(t)}$, $t\in\bbG$, where $q:\bbG\to L$ is the semigroup homomorphism defined earlier in Equation \eqref{eq:retract_onto_L}.

Since $\Alg[L]$ is biflat, the main results of this article imply that $L$ is uniformly locally finite, so that $(i)$ holds. Since $(i)$ holds, by Proposition~\ref{p:Schutz-is-iso} $\Alg[\bbG]$ is isomorphic as a Banach algebra to $\BG$; hence by Proposition~\ref{p:prelim_Cliff_biflat} every $G(x)$ is biflat, so that $(ii)$ holds as well.
\end{proof}

\subsection*{Outline proofs of the key propositions.}
\begin{proof}[Sketch proof of Proposition \ref{p:prelim_Cliff_biflat}]
Suppose $\BG$ is biflat. Let $x\in L$: then $\Alg[G(x)]$ is clearly a retract of $\BG$, and hence must itself be biflat. It is known that the $\lp{1}$-convolution algebra of a discrete group is biflat if and only if the group is amenable (see e.g.\ \cite[Thm~VII.2.35]{Hel_HBTA}), so $G(x)$ is amenable.

Suppose conversely that each group $G(x)$ is amenable. As remarked just above this implies that each $\Alg[G(x)]$ is biflat. In fact each $\Alg[G(x)]$ is $1$-biflat, i.e.~there exist contractive $\Alg[G(x)]$-bimodule maps $\sigma_x:\Alg[G(x)]\to (\Alg[G(x)]\ptp\Alg[G(x)])''$ such that $\pi''\sigma_x$ coincides with the embedding of $\Alg[G(x)]$ into its double dual. One then defines $\sigma: \BG \to (\BG\ptp \BG) ''$ to be the $\lp{1}$-sum (`coproduct') over $x$ of all the maps $\sigma_x$, and checks that it has the required properties as in Definition~\ref{dfn:biflat}.
\end{proof}

We turn to the task of proving Proposition \ref{p:Schutz-is-iso}. The required isomorphism is given by a more general notion of \Schutzrep\ which we outline briefly. For \emph{any partially ordered set $P$} define a bounded linear map $\Sch[P]:\lp{1}(P)\to \lp{\infty}(P)$ by
\[
 \Sch[P](e_t)={\bf 1}_{\dset{t}} \qquad\text{ for each $t\in P$}
\]
(cf.\ the definition in Section \ref{s:Schutz} for the special case $P=(L,\preceq)$).

Much of the formal calculation in previous sections survives in this more general setting. In particular, if $P$ is \emph{finite} then the formulas \eqref{eq:finitesum} and \eqref{eq:invSch} remain valid (with $L$ replaced by $P$); and the \emph{proof} of Theorem \ref{t:easydirection} goes through almost unchanged to give the following result.
\begin{propn}\label{p:gen_Sch_ULF}
Let $P$ be a uniformly locally finite, partially ordered set, i.e.\ one where
$\sup_{x\in P} \abs{\dset[P]{x}} < \infty$.
Then $\Sch[P]$ has range contained in $\lp{1}(P)$, and in fact is a continuous linear isomorphism from $\lp{1}(P)$ to $\lp{1}(P)$.
\end{propn}

To apply this proposition we introduce the following partial order on $\bbG$ (which is a special case of a more general construction for inverse semigroups). Given
$s,t\in \bbG$ we declare that
\[ s\preceq t \iff (\text{$s=tx$ for some $x \in L$})\;. \]
Note that $L$ is a downset in $(\bbG, \preceq)$ and that the partial order induced on $L$ coincides with its canonical partial order as a semilattice: thus the notation $(L,\preceq)$ is unambiguous.

Recalling that $t=tq(t)$ (see the earlier definition of the homomorphism $q:\bbG \to L$) we find that for every $t\in \bbG$
\[ \dset[\bbG]{t} = \{ tf \st f\in\dset[L]{q(t)} \} \;.\]
\emph{In particular, if the partially ordered set $(L,\preceq)$ is locally $C$-finite then so is $(\bbG,\preceq)$.}

\begin{proof}[Sketch proof of Proposition \ref{p:Schutz-is-iso}]
Suppose $L$ is uniformly locally finite. By the observation just made this implies $\bbG$ is uniformly locally finite; hence by Proposition \ref{p:gen_Sch_ULF} $\Sch[\bbG]$ is a bounded linear isomorphism from $\lp{1}(\bbG)$ onto itself. Now by the definitions of $\Alg[\bbG]$ and $\BG$ both have $\lp{1}(\bbG)$ as their underlying Banach space, so that we can consider $\Sch[\bbG]$ as a bounded linear isomorphism from $\Alg[\bbG]$ onto $\BG$.

\emph{Therefore it suffices to show that $\Sch[\bbG]$ is an algebra homomorphism.} 
We shall not give the details, but note that since $\Sch[\bbG]$ is known to be bounded linear it suffices to check the claim on the usual basis elements of $\Alg[\bbG]$.
\end{proof}

This completes our outline proofs of Propositions \ref{p:prelim_Cliff_biflat} and \ref{p:Schutz-is-iso}, and hence of Theorem \ref{t:Cliff_biflat}.

\subsection*{An aside on $\Sch[\bbG]$ for general Clifford semigroups}
The following remarks are not needed for our proof of Theorem \ref{t:Cliff_biflat} but provide extra background.

Recall that $\Sch[\bbG]:\Alg[\bbG] \to \lp{\infty}(\bbG)$ is a well-defined, bounded linear map for any Clifford semigroup $\bbG$, uniformly locally finite or otherwise. In general it does not seem possible to equip $\lp{\infty}(\bbG)$ with a Banach algebra structure \emph{with respect to which} $\Sch$ is an algebra homomorphism.

However, it turns out that $\Sch[\bbG]$ is in fact a bounded linear map taking values in the Banach space
\[ \lpsum{\infty}_{x\in \Und{L}} \lp{1}(G(x)) \]
and the same calculation as for the uniformly locally finite case will show that when this $\lp{\infty}$-sum is equipped with ``com\-ponent\-wise'' multiplication, $\Sch[\bbG]$ is an algebra homomorphism. This justifies speaking of the \emph{\Schutzrep}
\[ \Sch[\bbG]: \Alg[\bbG] \rTo \prod_{x\in \Und{L}} \Alg[G(x)] \]
for an \emph{arbitrary} Clifford semigroup $\bbG$, and generalises the case of finite Clifford semigroups which is presented in \cite{St_Mob}.

The proof that $\Sch[\bbG]$ takes values in $\lpsum{\infty}_{x\in \Und{L}} \lp{1}(G(x))$ is an easy consequence of the fact that for any given $t\in\bbG$ and $x\in L$, $\dset[\bbG]{t}$ intersects $G(x)$ in at most one point. This in turn follows from the following, more precise statement:

\medskip\noindent{\bf Claim.} Let $t\in \bbG$, $f\in L$. If $s\in \dset[\bbG]{t}\cap G(f)$, then $s=tf$.

\begin{proof}[Proof of claim]
Let $t\in \bbG$ and recall that $t\in G(q(t))$, $tq(t)=t$.

Suppose $s\in \dset[\bbG]{t}\cap G(f)$. Since $s\preceq t$ there exists $x\in L$ with $s=tx$; hence $s\in G(q(t))\cdot G(x)\subseteq G(q(t)x)$. Since $G(i)$ and $G(j)$ are disjoint for distinct $i,j\in L$ we must have $q(t)x=f$. Therefore $s=tx =tq(t)x = tf$ and the claim is proved.
\end{proof}


\appendix

\begin{section}{Proof of Proposition \ref{p:magic_idempotents}}

We follow the outline in \cite{Eti_Mob}, though our presentation is different. Recall that the idempotents $\rho_t$ are defined (for $t \in L$) by
\begin{equation}
\rho_t = \left\{ \begin{aligned}
	e_{\tht} & \quad\text{ if $t=\tht$} \\
	\prod_{x\in L \st x \prec t} (e_t-e_x) & \quad\text{ otherwise.}
 \end{aligned} \right.
\end{equation}
and that we wish to show that $\Sch(\rho_s)=\del_s$ for all $s \in L$. The case $s=\tht$ is trivial, so we shall henceforth assume that $s\succ \tht$.

If $c\in\Cplx^L$ and $s\in L$, we denote by $c_s$ the coefficient of $c$ with respect to the basis vector $\del_s$.

To prove that $\Sch(\rho_s)=\del_s$ it suffices to
 show that:
\begin{enumerate}
\item[(i)] $[\Sch(\rho_s)]_s=1$
\item[(ii)] $[\Sch(\rho_s)]_t=0$ for every $t \in L\setminus\{s\}$.
\end{enumerate}


We recall from the definition of $\Sch$ that given $x,y\in L$
\[ [\Sch(e_x)]_y = \left\{ \begin{aligned}
	1 & \quad\text{ if $x\leq y$} \\
	0 & \quad\text{ otherwise.}
 \end{aligned}\right. \]
and that since $\Sch$ is an algebra homomorphism,
 $[\Sch(a)\Sch(b)]_t=[\Sch(a)]_t[\Sch(b)]_t$ for all $t \in L$.

\begin{proof}[Proof of {\rm(i)}]
Expanding out the product that defines $\rho_s$ gives
\[ \rho_s = e_s + \sum_{x \in L \st x \prec s} \lambda_x e_x \]
for some scalar coefficients $\lambda_x$. Hence by the previous remarks
\[ [\Sch(\rho_s)]_s = [\Sch(e_s)]_s + \sum_{x \prec w} \lambda_x [\Sch(e_x)]_s = 1 \]
as required.
\end{proof}

\begin{proof}[Proof of {\rm(ii)}]
Let $t \in L\setminus\{s\}$. Then one of two cases must occur:

\begin{itemize}
\item {\bf Case A.} There exists $y \prec s$ with $t\preceq y$.
\item {\bf Case B.} For all $x \prec s$, $t \not\preceq x$.
\end{itemize}

Suppose we are in Case A. Then $\rho_s= b(e_s-e_y)$ for some $b \in \Alg[L]$, and so
\[ [\Sch(\rho_s)]_t = [\Sch(b)]_t\left([\Sch(e_s)]_t-[\Sch(e_y)]_t\right)
\;. \]
But since $t\preceq y\prec s$, $[\Sch(e_s)]_t = [\Sch(e_y)]_t=1$; therefore $[\Sch(\rho_s)]_t=0$.

On the other hand, suppose we are in Case B. We note that $\rho_se_s=\rho_s$.
Now since $st\preceq t$, the condition of Case B implies that $st \not\prec s$; hence $st=s$. Since we assumed that $s\neq t$ this forces $s\prec t$, and hence $[\Sch(e_s)]_t=0$. Therefore
\[ [\Sch(\rho_s)]_t = [\Sch(\rho_s)]_t[\Sch(e_s)]_t = 0 \;. \]

Thus for any $t \in L\setminus\{s\}$ we have $[\Sch(\rho_s)]_t=0$, and the proof of (ii) is complete.
\end{proof}
\end{section}

\subsection*{Acknowledgements}
This article grew out of work done while the author was a PhD student at the University of Newcastle upon Tyne. He would like to acknowledge the support of the School of Mathematics and Statistics while the article was written up, and to thank the referee for the close attention given and valuable suggestions made: in particular, as regards the present form of Theorem~\ref{t:easydirection} and Section~\ref{s:Clifford_sketch}.


\end{document}